\documentclass[12pt]{article}
\usepackage{amsmath, amssymb, amsthm}
\usepackage{mathrsfs}
\usepackage{amssymb}
\usepackage{amsmath}
\usepackage[hidelinks]{hyperref}
\usepackage{marvosym}
\usepackage{pict2e}
\usepackage{graphicx}
\usepackage{geometry}
\usepackage{color}

\def\q{\hfill\rule{1ex}{1ex}}
\def\0{\emptyset}

\def\p{{\bf Proof.}\quad}
\def\q{\hfill\rule{1ex}{1ex}}

\newtheorem*{transferlem}{Cocircuit-transfer lemma}

\newtheorem{theo}{Theorem}[section]
\newtheorem{coro}[theo]{Corollary}

\newtheorem{lem}[theo]{Lemma}

\makeatletter
\renewcommand{\@fnsymbol}[1]{\relax}

\makeatother

\begin{document}

\title{\bf Circuit spectra for binary matroids}

\author{Bing Wei~~~~~Haidong Wu~~~~~Shunzhe Zhang\\}

\author{Bing Wei\footnote{B. Wei: bwei@olemiss.edu}
~~~~~Haidong Wu\footnote{H. Wu: hwu@olemiss.edu}
~~~~~Shunzhe Zhang\footnote{S. Zhang: szhang12@olemiss.edu}\\
Department of Mathematics\\
The University of Mississippi\\
University, MS 38677, U.S.A\\
}

\date{}
\maketitle

\begin{abstract}
The circuit spectrum of a matroid is defined to be the set of circuit sizes. The circuit e-spectrum of a matroid is the set of circuit sizes of the matroid containing the element $e$. The circuit spectrum problem is an important problem. Erd\H{o}s (1992) conjectured that every connected graph $G$ with girth $g$ and average degree $d$ has a cycle spectrum of size at least $\Omega(d^{\lfloor (g-1)/2 \rfloor})$. Sudakov and Verstra\"{e}te~\cite{Su2008} proved the full conjecture of Erd\H{o}s in a stronger form.  It is natural to consider the size of the circuit spectra of binary matroids. Murty (1971) characterized all connected binary matroids with circuit spectra of size one. Lemos et al. (2011) characterized the connected binary matroids with circuit spectra of size two, where the largest size is odd. For general binary matroids, it is still an open question to find a lower bound for the size of the circuit spectra. In this paper, we give best-possible such lower bounds for connected binary matroids that do not contain both the Fano matroid and its dual as minors.  Let $M$ be a simple connected binary matroid that does not have both the Fano matroid and its dual as a minor, and $ d\geq 3$ be the minimum cocircuit size of $M$.   We prove that $(1)$ For each element of $M$, the circuit $e$-spectrum of $M$ must have size at least $d-2$. This bound is best possible. $(2)$ The circuit spectrum of $M$ must have size at least $d-1$, unless $M\cong F_7,F^*_7,R_{10}$. In exceptional cases, the spectrum of $M$ has size equal to $d-2$.

\vskip 0.1cm


\end{abstract}


\section{Introduction}

We use Oxley~\cite{Oxley1992} for terminology and notation not defined here. The {\it circuit spectrum} of a matroid $M$, denoted by $C(M)$, is the set of circuit sizes. The circuits of $M^*$ are called {\it cocircuits} of $M$. The {\it cocircuit spectrum} of a matroid $M$, denoted by $C^*(M)$, is the set of cocircuit sizes. We write $\mathcal{C}(M)$ for the set of circuits of $M$. So $C(M)=\{|C|:C\in\mathcal{C}(M)\}$ and $C^*(M)=C(M^*)$. For $e\in E(M)$, the circuit $e$-spectrum is $C_e(M)=\{|C|:C\in\mathcal{C}(M),\ e\in C\}$. A matroid is {\it binary} if and only if it has no minor isomorphic to the uniform matroid $U_{2,4}$. A matroid is {\it regular} if it is binary and has no minor isomorphic to the Fano matroid $F_7$, or its dual $F^*_7$.

Suppose that $M=M(G)$ is graphic. We use $C(G)$ to denote the {\it cycle spectrum} of $G$. For $S\subseteq E(G)$, the {\it cycle $e_S$-spectrum} of a graph $G$, denoted by $C_S(G)$, is the set of the sizes of cycles containing all edges in $S$. In particular, when $S=\{e\}$, we use the cycle $e$-spectrum to denote the cycle $e_{\{e\}}$-spectrum of $G$.

For a matroid $M$, the circuit spectrum of $M$ is the set of circuit sizes of $M$. Murty \cite{Mur1971} characterized all connected binary matroids with circuit spectra of size one. Lemos, Reid, and Wu \cite{Lrw2011} characterized the connected binary matroids with circuit spectra of size two, where the largest size is odd. Lewis, McNulty, Neudauer, Reid, and Sheppardson \cite{Reid2013} characterized the connected bicircular matroids with circuit spectra of size one. Bondy \cite{Bondy1971} showed that an $n$-vertex simple Hamiltonian graph with at least ${n^2\over 4}$ edges has cycles of every length unless it is isomorphic to $K_{\frac{n}{2},\frac{n}{2}}$. Motivated by this result, Beavers and Oxley \cite{Bo2005} showed that if $M$ is a simple rank-$r$ $GF(q)$-representable matroid with a very large number of elements, then the circuit spectrum of $M$ is $\{3, 4, \ldots, r+1\}$. In particular, a rank-$r$ simple binary matroid with at least $2^{r-1}$ elements either has circuits of all sizes or is isomorphic to $AG(r-1,2)$.

There is extensive literature on the cycle spectra of graphs (see, for example, \cite{Er1999}, \cite{Liu2018}). Erd\H{o}s (1992) conjectured that every connected graph $G$ with girth $g$ and average degree $d$ has a cycle spectrum of size at least $\Omega(d^{\lfloor (g-1)/2 \rfloor})$. Sudakov and Verstra\"{e}te~\cite{Su2008} proved the full conjecture of Erd\H{o}s in a stronger form.

Suppose that $M=M^*(G)$ is cographic where $G$ is simple with girth $g$. Then the circuit spectrum of $M$ is the {\it bond spectrum} of $G$, which is the set of sizes of bonds (minimal non-empty edge-cuts) of $G$. We denote this spectrum by $C^*(G)$. Thus $C^*(G)=C(M^*(G))$, whereas $C(G)=C(M(G))$.



In a previous paper \cite{Wei2025}, we gave an exponential lower bound on the bond spectra for connected graphs with minimum degree at least three and with girths at least seven. In \cite{Wei2025-2}, we gave the best possible lower bound on the bond spectra for connected graphs with minimum degree at least three, and we characterized all connected graphs attaining the lower bound. For binary matroids, a general lower bound for the size of the circuit spectrum is not known. Several years ago, the second author of the current paper asked the following question: 

\noindent {\bf Question:} Let $M$ be a simple connected binary matroid with minimum cocircuit size $d\geq 3$. Find a function $f(d)$ such that the size of the circuit spectrum of $M$ is at least $f(d)$.  


Let $G$ be a graph. A set of $k$ paths $P_1, P_2,\ldots, P_k$ are {\it admissible} if $|E(P_1)|\geq 2$ and $|E(P_1)|, |E(P_2)|,\ldots, |E(P_k)|$ form an arithmetic progression of length $k$ with common difference one or two.  J. Gao, Q. Huo, C. Liu and J. Ma \cite{Gao2022} proved the following result for graphic matroids.

\begin{theo}\label{1.1} {\rm (\cite{Gao2022})}
Let $G$ be a 2-connected graph and let $x, y$ be distinct vertices of $G$. If every vertex of $G$ other than $x$ and $y$ has degree at least $k+1$, then there exist $k$ admissible paths from $x$ to $y$ in $G$, where $k$ is a positive integer.
\end{theo}

As a direct consequence, for any edge $e$ in a 2-connected graph with minimum degree $d$, the $e$-spectrum of $G$ contains a sequence of at least $d-1$ circuits so that the consecutive terms have size differences by one or two.



Gao et al. \cite{Gao2022} showed that the above result is best-possible. For non-bipartite graphs, consider the complete graph $K_{d+1}$.  Let $e$ be an edge of $K_{d+1}$. Then $K_{d+1}$ satisfies the minimum degree condition and has exactly $d-1$ distinct cycle lengths containing $e$ with consecutive lengths. For bipartite graphs, consider the complete bipartite graph $K_{d,n}$ for $n\geq d$ and $x,y$ are two distinct vertices in the partite set of size $d$. Then $K_{d,n}$ satisfies the minimum degree condition, and the lengths of $x$-$y$ paths in $K_{d,n}$ take exactly $d-1$ values, namely $2,4,\ldots,2d-2$. Both examples show that the bounds of Theorem \ref{1.1} are tight. If, instead of requiring the circuits to have size differences by one or two, we only require that the circuits have different sizes, then the above examples also show that $d-1$ is the best possible lower bound.

Murty \cite{Mur1971} characterized all connected binary matroids with circuit spectra of size one.

\begin{theo}\label{Murty} {\rm (Murty~\cite{Mur1971})} 
Let $M$ be a connected binary matroid. For a positive integer $c$, $C(M)=\{c\}$ if and only if $M$ is isomorphic to one of the following matroids:

$(i)$ a $c$-subdivision of $U_{0,1}$,

$(ii)$ a $k$-subdivision of $U_{1,n}$, where $c=2k$ and $n\geq 3$,

$(iii)$ an $l$-subdivision of $PG(r,2)^*$, where $c=2^rl$ and $r\geq 2$,

$(iv)$ an $l$-subdivision of $AG(r+1,2)^*$, where $c=2^rl$ and $r\geq 2$.
\end{theo}

Lemos, Reid, and Wu \cite{Lrw2011} characterized the connected binary matroids with circuit spectra of size two, where the largest size is odd.

\begin{theo}\label{LRW} {\rm (Lemos, Reid, and Wu~\cite{Lrw2011})} 
Let $c,d\in\mathbb{Z}^{+}$ with $c<d$ and $d$ odd.
Let $M$ be a connected binary matroid. Then
$C(M)=\{c,d\}$ if and only if there are connected
binary matroids $M_0,M_1,\ldots,M_n$ for some $n\in\mathbb{Z}^{+}$
such that the following hold.

$(i)$ $E(M_i)\cap E(M_j)=\{e\}$, for distinct $i$ and $j$ in $\{0,1,\ldots,n\}$.

$(ii)$ $E(M_0)$ is a circuit of $M_0$.

$(iii)$ For $i\in\{1,2,\ldots,n\}$, $\{e\}$ is a series class of $M_i$, all other series classes of $M_i$ have size $l_i$, and the cosimplification of $M_i$ is isomorphic to one of the following matroids.
     
    $(a)$ $U_{1,n_i}$, for some $n_i\geq 3$, where $c=2l_i$.

    $(b)$ $PG(r_i,2)^*$, for some $r_i\geq 2$, where $c=2^{r_i}l_i$.

    $(c)$ $AG(r_i,2)^*$, for some $r_i\geq 3$, where $c=2^{r_i-1}l_i$.

    $(d)$ $S_{n_i}^*$, for some $n_i\geq 4$, and $e$ is the cotip, where $c=4l_i$.

    $(e)$ $B(r_i,2)^*$, for some $r_i\geq 3$, and $e$ is the cotip, where $c=2^{r_i}l_i$.

$(iv)$ $d=|E(M_0)|-1+d_1+d_2+\cdots+d_n>c$, where $d_i=\frac{c}{2}$ when $(iii)(a)$ holds, $d_i=(2^{r_i}-1)l_i$ when $(iii)(b)$ holds, $d_i=(2^{r_i-1}-1)l_i$ when $(iii)(c)$ holds, $d_i=n_il_i$ when $(iii)(d)$ holds, and $d_i=c$ when $(iii)(e)$ holds.

$(v)$  $M=S(M_0,M_1,\ldots,M_n)/e$.
\end{theo}

In \cite{Lrw2010}, Lemos, Reid, and  Wu proved the following result. 

\begin{theo}\label{e-spectrum-matroid} 
Suppose that $e$ is an element of a binary matroid $M$. If no $F_7^*$ -minor of $M$ uses $e$, then the following statements are equivalent:

\noindent (i) There is a subdivision $N$ of $M$ such that the circuit $e$-spectrum of $N$ has size one, and each circuit of $N$ containing $e$ has size equal to the size of the largest circuit of $M$ containing $e$. \\
(ii) $M$ is a series-parallel network. 
\end{theo}

P. Li, L. Wang, Y. Wu, and H.J. Lai \cite{Lai2019} gave a characterization of all connected binary rooted matroids with an element lying in at most three circuits, and a characterization of all connected binary matroids such that some element lies in all but at most three circuits.

The above theorems show that there are infinitely many binary matroids whose circuit spectra have size at most two.  Moreover, $AG(2k+1,2)$ has precisely  $k$ circuit sizes $4,6,\ldots,2k+2 $ and cogirth  $d=4^k$. For example, $AG(5,2)$ has circuit sizes $\{4,6\}$ and cogirth $16$. Therefore, we do not expect the circuit spectrum of a general binary matroid to be a linear function of the cogirth $d$.

It is natural to ask what lower bounds can be obtained for the size of the circuit spectra of binary matroids that do not have both the Fano matroid $F_7$ and its dual as a minor. In this paper, we completely resolve this problem. We prove a linear lower bound of the minimum cocircuit size $d$ on the size of the circuit spectra for simple and cosimple connected binary matroids that do not have both the Fano matroid $F_7$ and its dual as a minor. We also give a best possible lower bound for the $e$-spectrum for such binary matroids.

In our first result, we provide the bond $e$-spectra of simple 2- and 3-connected graphs, that is, the $e$-spectra of connected and 3-connected cographic matroids.

\begin{theo}\label{1.2} 
$(a)$ Let $G$ be a simple 2-connected graph with minimum degree $\delta\geq 3$ and $e\in E(G)$. Suppose that $g(G-e)=d\geq 3$ and $h=d+\lfloor \frac{d-1}{2} \rfloor-2$. Then the size of the bond $e$-spectrum of $G$ is at least $h$.

$(b)$ Let $G$ be a simple 3-connected graph and $T=\{e,f,g\}$ be a 3-cocircuit in $G$. Suppose that $g(G-T)\geq d\geq 3$ and $h=d+\lfloor \frac{d-1}{2} \rfloor-2$. If $G \ncong K_4$, then there exist bonds $C_i$ and $D_j$ of $G$ for $i\in \{1,\ldots, h\}$, and $j\in \{1,\ldots, h-1\}$ such that $C_i\cap T=\{e\}$ and $D_j\cap T=\{f,g\}$, with $|C_1|,\ldots,|C_{h}|$ pairwise distinct and $|D_1|,\ldots,|D_{h-1}|$ pairwise distinct.
\end{theo}

The bounds in Theorem \ref{1.2} are attained in some cases, including the complete bipartite graph $K_{3,3}$ and the Petersen graph $P_{5,2}$. 


The following two theorems are the main results of this paper. 

\begin{theo}\label{1.3} 
Let $M$ be a simple connected binary matroid that does not have both an $F_7$-minor and an $F_7^*$-minor, and let $e$ be an element of $M$ with $|E(M)|\geq 2$. Suppose that every cocircuit $Y$ of $M$ with $e\notin Y$ has size at least $d\geq 3$. Then the size of the circuit $e$-spectrum of $M$ is at least $d-2$.
\end{theo}

It is easy to show that when $M\cong R_{10}$, the lower bound $d-2$ is attained in Theorem \ref{1.3}. Take $k$ copies of the Fano dual $F_7^*$ with a fixed element $e$. Then take the parallel connection of these $k$ copies along the common element $e$. The resulting matroid, denoted by $P_k(F_7^*, e)$, is clearly connected with minimum cocircuit size $d=3$. It also has no $F_7$-minor. It is straightforward to see that the circuit $e$-spectrum of the element $e$ has size $d-2=1$. This shows that the lower bound in Theorem \ref{1.3} is sharp. 

In the next theorem, we give a best possible bound on the sizes of the spectra of binary matroids that do not have both an $F_7$-minor and an $F_7^*$-minor.

\begin{theo}\label{1.4} 
Let $M$ be a nonempty simple connected binary matroid that does not have both an $F_7$-minor and an $F_7^*$-minor. Suppose that every cocircuit of $M$ has size at least $d\geq 3$. Then the size of the circuit spectrum of $M$ is at least $d-1$, unless $M\cong F_7,F^*_7,R_{10}$. In the exceptional cases, the size of the circuit spectrum of $M$ is exactly $d-2$ if $d$ is the minimum cocircuit size of $M$.
\end{theo}

Let $M_1, M_2$ be binary matroids with element sets $S_1,S_2$ respectively, where $S_1,S_2$ may intersect. Seymour~\cite{Se1980} defined a new binary matroid $M_1\Delta M_2$ to be the matroid with element set $S_1\Delta S_2$ and with cycles all subsets of $S_1\Delta S_2$ of the form $C_1\Delta C_2$, where $C_i$ is a cycle of $M_i$ for $i\in \{1,2\}$. (For sets $S_1,S_2$, $S_1\Delta S_2$ denotes $(S_1-S_2)\cup (S_2-S_1)$. A cycle of a binary matroid is a subset of the elements expressible as a disjoint union of circuits. It is easy to see that if $C, C'$ are cycles then $C\Delta C'$ is a cycle.) Seymour~\cite{Se1980} was concerned with three special cases of this operation, as follows. 

$(i)$ When $S_1\cap S_2=\emptyset$ and $S_1,S_2\neq \emptyset$. $M_1\Delta M_2$ is a 1-sum of $M_1$ and $M_2$.

$(ii)$ When $|S_1\cap S_2|=1$, and $S_1\cap S_2=\{z\}$ and $z$ is not a loop or 
coloop of $M_1$ or $M_2$, and $|S_1|,|S_2|\geq 3$, $M_1\Delta M_2$ is a 2-sum of $M_1$ and $M_2$.

$(iii)$ When $|S_1\cap S_2|=3$, and $S_1\cap S_2=Z$ and $Z$ is a circuit of $M_1$ and $M_2$, and $Z$ includes no cocircuit of either $M_1$ or $M_2$, and $|S_1|,|S_2|\geq 7$, $M_1\Delta M_2$ is a 3-sum of $M_1$ and $M_2$.

If $(X_1, X_2)$ is a partition of $E(M)$, we say that $(X_1, X_2)$ is a {\it $k$-separation} of $M$ if $|X_1|,|X_2|\geq k$ and $r(X_1)+r(X_2)\leq r(M)+k-1$. A $k$-separation $(X_1, X_2)$ is {\it exact} if $|X_1|,|X_2|\geq k$ and $r(X_1)+r(X_2)=r(M)+k-1$. $M$ is {\it $s$-connected} if, for all $k\in \{1,\ldots,s-1\}$, it has no $k$-separations. A {\it $k$-subdivision} of a matroid is obtained by replacing each element by a series class of size $k$.

\section{Some Preliminaries}

In this section, we list some theorems and prove preliminary results which will be used in our proofs.

\begin{theo}\label{2.3} (Seymour~\cite[(2.6)]{Se1980})
If $(X_1,X_2)$ is an exact 2-separation of $M$, then there are matroids 
$M_1,M_2$ on $X_1\cup \{z\}$, $X_2\cup \{z\}$, respectively (where $z$ is a new element), such that $M$ is the 2-sum of $M_1$ and $M_2$. Conversely, if $M$ is the 2-sum of $M_1$ and $M_2$, then $(E(M_1)-E(M_2),E(M_2)-E(M_1))$  is an exact 2-separation of $M$, and $M_1,M_2$ are isomorphic to minors of $M$.
\end{theo}

The following theorem was proved by Seymour~\cite{Se1980}. We use the following restatement in Oxley~\cite{Oxley1992}. 

\begin{theo}\label{2.4} 
Let $(X_1,X_2)$ be an exact 3-separation of a binary matroid $M$ with $|X_1|$, $|X_2|\geq 4$, and let $Z$ be a 3-element set that is disjoint from $E(M)$. Then there is a binary matroid $N$ with ground set $E(M)\cup Z$ such that $N|E(M)=M$, the set $Z$ is a circuit of $N$, and $r(X_i\cup Z)=r(X_i)$ for each $i$. Moreover, $M=N_1\oplus_3 N_2$ where $N_i=N|(X_i\cup Z)$ for each $i$, and $N$ is unique up to a permutation on the labels of the elements of $Z$. Conversely, the 3-sum, $M_1\oplus_3 M_2$, of binary matroids $M_1$ and $M_2$ has $(E(M_1)-E(M_2),E(M_2)-E(M_1))$ as an exact 3-separation, $|E(M_1)-E(M_2)|$ and $|E(M_2)-E(M_1)|$ exceed three, and
$r(M_1\oplus_3 M_2)=r(M_1)+r(M_2)-2$. Moreover, if $M$ is 3-connected, then $N_1$ and $N_2$ are isomorphic to minors of $M$ \cite{Se1980}.
\end{theo}


\begin{theo}\label{2.5} {\rm (Seymour~\cite[(7.6)]{Se1980})}
$F_7$ is a splitter for the class of binary matroids without $F^*_7$ minors, and dually, $F_7^*$ is a splitter for the class of binary matroids without $F_7$ minors.
\end{theo}

\begin{theo}\label{2.6} {\rm (Seymour~\cite[(14.3)]{Se1980})}
Every regular matroid $M$ may be constructed by means of 1-sum, 2-sums, 
and 3-sums, starting with matroids each isomorphic to a minor of $M$ and 
each either graphic or cographic or isomorphic to $R_{10}$.
\end{theo}

For sets $A,B\subseteq\mathbb{Z}$, their sumset is
$A+B=\{a+b:a\in A,\ b\in B\}$. We will use the following elementary inequality.

\begin{lem}\label{Cauchy}
Let $A$ and $B$ be finite nonempty subsets of $\mathbb{Z}$. Then $|A+B|\geq |A|+|B|-1.$ \end{lem}

\p Write $A=\{a_1<a_2<\cdots<a_m\}$ and
$B=\{b_1<b_2<\cdots<b_n\}$. The strictly increasing sequence $a_1+b_1<a_2+b_1<\cdots<a_m+b_1<a_m+b_2<\cdots<a_m+b_n$ contains $m+n-1=|A|+|B|-1$ distinct elements of $A+B$.
The desired inequality follows. \q

We denote by $E[X,Y]$ the set of edges of $G$ with one end in $X$ and the other end in $Y$. A subset $Y$ of $V(G)$ is denoted by $\overline{X}$ if $Y=V(G)\setminus X$. The set $E[X,\overline{X}]$ is called the \emph{edge-cut} of $G$ associated with $X$ and denoted by $\partial(X)$. A subgraph $H$ of $G$ is called \emph{non-separating} if $G-H$ is connected.

\begin{lem}\label{2.7}
Let $G$ be a connected graph with minimum degree at least three and $e=uv\in E(G)$. Suppose that $G$ has a non-separating induced tree $T$ with $t$ vertices such that $u\in V(T)$ and $v\notin V(T)$. Then the size of the bond $e$-spectrum of $G$ is at least $t$.
\end{lem}

\p Let $T=T_t$ be such a tree with $V(T_t)=\{x_1,x_2,\ldots,x_t\}$. Assume that $x_1=u$ and $T_{i}=T_{i+1}-\{x_{i+1}\}$ for $i\in \{1,2,\ldots,t-1\}$, where $x_{i+1}$ is a leaf of $T_{i+1}$ other than $x_1$. Then $T_i$ is an induced tree of $G$ with $i$ vertices such that $x_1\in V(T_i)$ and $v\notin V(T_i)$. Moreover, $T_i$ is non-separating: $G-V(T_{i+1})$ is connected, and since $x_{i+1}$ has exactly one neighbor in $T_{i+1}$ and degree at least three, it has a neighbor in $G-V(T_{i+1})$. Hence $\partial(T_i)$ is a bond of $G$ with $e\in \partial(T_i)$. Since $T_i$ is an induced tree, $|\partial(T_i)|=\sum_{s=1}^{i}d_G(x_s)-2|E(T_i)| =\sum_{s=1}^{i}d_G(x_s)-2i+2.$ Since $\delta(G)\geq 3$, we have  $2<|\partial(T_{i})|<|\partial(T_{i+1})|$ for $i\in \{1,2,\ldots,t-1\}$. This implies $\{|\partial(T_{i})|: i\in \{1,\ldots, t\}\}\subseteq C^*(G)$ and thus the size of the bond $e$-spectrum of $G$ is at least $t$.
\q

\section{Graphic and Cographic Matroids}

The following are two consequences of Theorem \ref{1.1} of J. Gao, Q. Huo, C. Liu, and J. Ma~\cite{Gao2022}.

\begin{coro}\label{3.1} 
Let $G$ be a 2-connected graph and $e=uv\in E(G)$. If every vertex in $G$ other than $u,v$ has degree at least $d\geq 3$, then the size of the cycle $e$-spectrum of $G$ is at least $d-1$. 
\end{coro}

\begin{coro}\label{3.2} 
Let $G$ be a 3-connected graph and $T=uevfwgu$ be a triangle in $G$, where $e=uv$, $f=vw$ and $g=wu$ are edges of $G$. If every vertex in $G$ other than $u,v,w$ has degree at least $d\geq 3$, then there exist cycles $C_i,D_i$ of $G$ for $i\in \{1,\ldots, d-2\}$ such that $C_i\cap T=\{e\}$ and $D_i\cap T=\{f,g\}$, with $|C_1|,\ldots,|C_{d-2}|$ pairwise distinct and $|D_1|,\ldots,|D_{d-2}|$ pairwise distinct.
\end{coro}

\p Let $H=G-w$. Then $H$ is 2-connected and every vertex in $H$ other than $u,v$ has degree at least $d-1\geq 2$. By Theorem \ref{1.1}, $H$ contains $d-2$ admissible paths $P_i$ from $u$ to $v$, where $i\in \{1,\ldots,d-2\}$. Note that $|E(P_i)|\geq 2$. Set $C_i=P_i\cup \{e\}$ and $D_i=P_i\cup \{f,g\}$. Thus $C_i\cap T=\{e\}$ and $D_i\cap T=\{f,g\}$. 
\q


Now, we will prove Theorem \ref{1.2}. Before we prove Theorem \ref{1.2}, we prove the next results.

\begin{theo}\label{3.3} 
Let $G$ be a simple 2-connected graph with minimum degree $\delta\geq 3$ and $e=uv\in E(G)$. Suppose that $g(G-e)=d\geq 3$ and $h=d+\lfloor \frac{d-1}{2} \rfloor-2$. Then $G$ has a non-separating induced tree $T$ with $|V(T)|\geq h$ such that $u\in V(T)$ and $v\notin V(T)$. Furthermore, the size of the bond $e$-spectrum of $G$ is at least $h$.
\end{theo}

\p Since $g(G-e)=d$ and $e=uv$, every circuit of $G-e$ has length at least $d$. Since $\lfloor \frac{d-1}{2} \rfloor=\lceil \frac{d}{2} \rceil-1\leq \frac{d-1}{2}$, we have $2h\leq 3d-5$, and $h-\lceil \tfrac{d}{2} \rceil=d-3$.

We call an induced tree $T$ of $G$ {\it good} if $u\in V(T)$, $v\notin V(T)$ and $G-V(T)$ is connected. Since $G-u$ is connected, $G[\{u\}]$ is good, and we may choose a good tree $T$ with $|V(T)|$ as large as possible. Suppose, to the contrary, that $|V(T)|\leq h-1$. So $|E(T)|\leq h-2$. Let $H=G-V(T)$. For a graph $J$, let $X_J$ be the set of cut vertices of $J$.

\noindent\textbf{Claim A.} Let $U\subseteq V(G)\setminus\{v\}$ be non-empty such that $J=G-U$ is connected. If no vertex of $V(J)\setminus (\{v\}\cup X_J)$ has a neighbor in $U$, then $V(J)=\{v\}$.

\p Suppose that $|V(J)|\geq 2$. If $X_J=\emptyset$, then some vertex of $V(J)\setminus\{v\}$ has a neighbor in $U$ as $G-v$ is connected, a contradiction. Thus $X_J\neq \emptyset$ and $J$ has at least two end blocks. We choose an end block $B$ with a cut vertex $w$ such that $v\notin V(B)\setminus\{w\}$. Since $G-w$ is connected, one vertex of $V(B)\setminus\{w\}$ has a neighbor in $U$, a contradiction. \q

If $V(H)=\{v\}$, then at least $3|V(T)|-2|E(T)|=|V(T)|+2$ edges join $V(T)$ to $v$ as $\delta\geq 3$, a contradiction. Hence $V(H)\neq\{v\}$. By Claim A, applied with $U=V(T)$ and $J=H$, there is a vertex $x_1\in V(H)\setminus (\{v\}\cup X_H)$ with a neighbor in $T$. Let $Q=V(T)\cup\{x_1\}$ and $H_0=H-x_1$. 

\noindent\textbf{Claim B.} $x_1$ has exactly two neighbors $y_1,z_1$ in $T$. Consequently, $|Q|\leq h$, $G[Q]$ is connected and has unique circuit $C_1=x_1y_1Tz_1x_1$.

\p If $x_1$ has exactly one neighbor in $T$, then $G[V(T)\cup\{x_1\}]$ is a good tree as $H-x_1$ is connected, which contradicts the maximality of $|V(T)|$. If $x_1$ has three neighbors $y,z,w$ in $T$, then $y,z,w\in V(T)$, and $x_1yTzx_1$, $x_1zTwx_1$ and $x_1yTwx_1$ are circuits of $G-v$ of lengths $|E(yTz)|+2$, $|E(zTw)|+2$ and $|E(yTw)|+2$, respectively. Since every edge of the tree $T$ lies on at most two of the paths $yTz$, $zTw$ and $yTw$, the sum of their lengths is at most $2|E(T)|+6\leq 2h+2\leq 3d-3$, which contradicts the fact that each of them has length at least $d$. Hence $x_1$ has exactly two neighbors $y_1,z_1$ in $T$, and $y_1,z_1\in V(T)$. Since $T$ is a tree and $H-x_1$ is connected, we have $|Q|=|V(T)|+1\leq h$, $G[Q]$ is connected and has a unique circuit $C_1=x_1y_1Tz_1x_1$. \q

\noindent\textbf{Claim C.} For every $w\in V(H_0)$, at most one edge of $G-e$ joins $w$ to $Q$.

\p By Claim B, $G[Q]-E(C_1)$ is a forest each of whose components contains exactly one vertex of $C_1$. Let $a,b\in Q$, and let $a',b'$ be the vertices of $C_1$ in the components containing $a$ and $b$, respectively. The paths from $a$ to $a'$ and from $b'$ to $b$ in this forest, together with a shortest path from $a'$ to $b'$ in $C_1$, contain a path from $a$ to $b$ with at most $|Q\setminus V(C_1)|+\lfloor \tfrac{|V(C_1)|}{2} \rfloor+1=|Q|-\lceil \tfrac{|V(C_1)|}{2} \rceil+1\leq h-\lceil \tfrac{d}{2} \rceil+1=d-2$ vertices. If two edges $wa,wb$ of $G-e$ joined $w$ to $Q$, then $w$ with such an $(a,b)$-path would form a circuit of $G-e$ with at most $d-1$ vertices, which contradicts $g(G-e)=d$. \q

If $V(H_0)=\{v\}$, then by Claim B, $G[Q]$ has exactly $|Q|$ vertices and $|Q|$ edges. Since $\delta\geq 3$, there are at least $3|Q|-2|Q|=|Q|\geq |V(C_1)|\geq 3$ edges joining $Q$ to $v$, which contradicts Claim C. Hence $V(H_0)\neq\{v\}$. By Claim A, applied with $U=Q$ and $J=H_0$, there is a vertex $w\in V(H_0)\setminus (\{v\}\cup X_{H_0})$ with a neighbor in $Q$. By Claim C, let $N(w)\cap Q=\{a\}$ for some $a$. By Claim B, $G[Q]-E(C_1)$ is a forest each of whose components contains exactly one vertex of $C_1$. Hence at most $|Q\setminus V(C_1)|\leq h-d$ vertices of $C_1$ have degree at least three in $G[Q]$. Since $|V(C_1)|-(h-d)\geq 2d-h=d-\lfloor \frac{d-1}{2} \rfloor+2\geq 3$, there is a vertex $z\in V(C_1)\setminus\{a,u\}$ with $d_{G[Q]}(z)=2$. Thus $G[Q\setminus\{z\}]$ is a tree containing $u$, and $T^*=G[(Q\setminus\{z\})\cup\{w\}]$ is an induced tree with $u\in V(T^*)$ and $v\notin V(T^*)$. Since $\delta\geq 3$ and Claim C, $z$ has a neighbor in $H_0$, which is different from $w$. Since $H_0-w$ is connected, we have $G-V(T^*)=G[(V(H_0)\setminus\{w\})\cup\{z\}]$ is also connected. Hence $T^*$ is a good tree with $|V(T^*)|=|V(T)|+1$, which contradicts the maximality of $|V(T)|$.

Therefore $G$ has a non-separating induced tree $T$ with $|V(T)|\geq h$ such that $u\in V(T)$ and $v\notin V(T)$. By Lemma \ref{2.7}, the size of the bond $e$-spectrum of $G$ is at least $h$. This completes the proof of Theorem \ref{3.3}.
\q

The next two lemmas, Lemmas~\ref{3.4} and~\ref{3.5}, generalize Lemma~8(a) and Lemma~8(b) of Hochst\"{a}ttler and Jackson~\cite{Hoch1998}, respectively.

\begin{lem}\label{3.4}
Let $G$ be a simple 3-connected graph with minimum degree $\delta=3$ and $v\in V(G)$ with $N(v)=\{v_1,v_2,v_3\}$. Suppose that $g(G-v)=d\geq 3$ and $h=d+\lfloor \frac{d-1}{2} \rfloor-2$. Then $G$ has a non-separating induced tree $T$ with $|V(T)|\geq h$ such that $v,v_3\in V(T)$ and $v_1,v_2\notin V(T)$.
\end{lem}

\p If $d=3$, then $h=2$, and $G[\{v,v_3\}]$ is the required tree as $G-\{v,v_3\}$ is connected. Thus we may assume that $d\geq 4$. Since $g(G-v)=d$, every circuit of $G-v$ has length at least $d$. Since $\lfloor \frac{d-1}{2} \rfloor=\lceil \frac{d}{2} \rceil-1\leq \frac{d-1}{2}$, we have $2h\leq 3d-5$, and $h-\lceil \tfrac{d}{2} \rceil=d-3$.

Let $Z=\{v_1,v_2\}$. We call an induced tree $T$ of $G$ {\it good} if $v,v_3\in V(T)$, $V(T)\cap Z=\emptyset$ and $G-V(T)$ is connected. Since $G[\{v,v_3\}]$ is good, we may choose a good tree $T$ with $|V(T)|$ as large as possible. Suppose, to the contrary, that $|V(T)|\leq h-1$. Let $T'=T-v$ and $H=G-V(T)$. Since $N(v)=\{v_1,v_2,v_3\}$, the vertex $v$ is a leaf of $T$ whose only neighbor in $T$ is $v_3$, and $Z$ is the set of vertices of $H$ adjacent to $v$. Thus $T'$ is a tree with $|E(T')|=|V(T)|-2\leq h-3$. For a graph $J$, let $X_J$ be the set of cut vertices of $J$.

\noindent\textbf{Claim A.} Let $U\subseteq V(G)\setminus Z$ be non-empty such that $J=G-U$ is connected. If no vertex of $V(J)\setminus (Z\cup X_J)$ has a neighbor in $U$, then $J$ is a path with ends $v_1$ and $v_2$.

\p Let $B$ be a block of $J$, let $X_B$ be the set of cut vertices of $J$ that lie in $B$, and let $P_B$ be the set of vertices of $B$ which are not cut vertices of $J$. Every vertex of $P_B$ has all its neighbors in $V(B)\cup U$, and the sets $P_B$ are pairwise disjoint for distinct blocks $B$. Suppose that $S=P_B\setminus Z\neq\emptyset$. By hypothesis, no vertex of $S$ has a neighbor in $U$. So $N_G(S)\subseteq X_B\cup (Z\cap V(B))$. Since $N_G(S)$ separates $S$ from $U$ in $G$ and $G$ is 3-connected, we obtain the following.

$(\ast)$ If $P_B\not\subseteq Z$, then $|X_B\cup (Z\cap V(B))|\geq 3$.

If $X_J=\emptyset$, then $J$ is a block with $P_J=V(J)$. So $|X_B\cup (Z\cap V(B))|\leq 2$. By $(\ast)$, we get $V(J)\subseteq Z$ and $J$ is the edge $v_1v_2$, so the claim holds. Otherwise, $X_J\neq \emptyset$ and $J$ has at least two end blocks. Each end block $B$ contains exactly one cut vertex of $J$, and thus $P_B\neq\emptyset$. If one of end blocks satisfies $P_B\cap Z=\emptyset$, then $|X_B\cup (Z\cap V(B))|=1$, which contradicts $(\ast)$. Thus each end block satisfies $P_B\cap Z\neq\emptyset$. Note that $|Z|=2$. Then $J$ has exactly two end blocks $B_1,B_2$, and we may assume that $P_{B_i}\cap Z=\{v_i\}$ for $i\in\{1,2\}$. This implies that $|X_{B_i}\cup (Z\cap V(B_i))|\leq 2$. By $(\ast)$, we get $P_{B_i}=\{v_i\}$, and $B_i$ is an edge. Note that every other block $B$ of $J$ contains exactly two cut vertices of $J$ and $Z\cap V(B)=\emptyset$. By $(\ast)$, we get $P_B=\emptyset$, and thus $B$ is an edge. Therefore $J$ is a path with ends $v_1$ and $v_2$. \q

\noindent\textbf{Claim B.} $H$ is not a path with ends $v_1$ and $v_2$.

\p Suppose that $H$ is a path with ends $v_1$ and $v_2$. Since $\delta=3$, every vertex of $H$ has a neighbor in $T'$. If $|V(H)|=2$, then $T$ has a leaf $u\neq v$ with at least two neighbors outside $T$. So $uv_1v_2u$ is a triangle of $G-v$, which contradicts $d\geq 4$. Thus $H$ has three consecutive vertices $c_1,c_2,c_3$. Choose $u_i\in N(c_i)\cap V(T')$ for $i\in\{1,2,3\}$. Then $c_1c_2u_2T'u_1c_1$, $c_2c_3u_3T'u_2c_2$ and $c_1c_2c_3u_3T'u_1c_1$ are circuits of $G-v$ (where $u_iT'u_j$ is trivial if $u_i=u_j$) of lengths $|E(u_1T'u_2)|+3$, $|E(u_2T'u_3)|+3$ and $|E(u_1T'u_3)|+4$, respectively. Since every edge of $T'$ lies on at most two of the paths $u_1T'u_2$, $u_2T'u_3$ and $u_1T'u_3$, the sum of these three lengths is at most $2|E(T')|+10\leq 2h+4\leq 3d-1$. This contradicts $g(G-v)=d\geq 4$. \q

By Claims A and B, applied with $U=V(T)$ and $J=H$, there is a vertex $x_1\in V(H)\setminus (Z\cup X_H)$ with a neighbor in $T$. Let $Q=V(T)\cup\{x_1\}$, $Q^\circ=Q\setminus\{v\}$, and $H_0=H-x_1$. 

\noindent\textbf{Claim C.} $x_1$ has exactly two neighbors $y_1,z_1$ in $T$. Consequently, $|Q|\leq h$, $G[Q]$ and $G[Q^\circ]$ are connected and have a unique circuit $C_1=x_1y_1T'z_1x_1$.

\p If $x_1$ has exactly one neighbor in $T$, then $G[V(T)\cup\{x_1\}]$ is a good tree as $H-x_1$ is connected, which contradicts the maximality of $|V(T)|$. If $x_1$ has three neighbors $y,z,w$ in $T$, then $y,z,w\in V(T')$, and $x_1yT'zx_1$, $x_1zT'wx_1$ and $x_1yT'wx_1$ are circuits of $G-v$ of lengths $|E(yT'z)|+2$, $|E(zT'w)|+2$ and $|E(yT'w)|+2$, respectively. Since every edge of the tree $T'$ lies on at most two of the paths $yT'z$, $zT'w$ and $yT'w$, the sum of their lengths is at most $2|E(T')|+6\leq 2h\leq 3d-5$, which contradicts the fact that each of them has length at least $d$. Hence $x_1$ has exactly two neighbors $y_1,z_1$ in $T$, and $y_1,z_1\in V(T')$. Since $T$ is a tree and $H-x_1$ is connected, we have $|Q|=|V(T)|+1\leq h$, $G[Q]$ and $G[Q^\circ]$ are connected and has a unique circuit $C_1=x_1y_1T'z_1x_1$. \q

\noindent\textbf{Claim D.} Any two vertices of $Q^\circ$ are joined by a path of $G[Q^\circ]$ with at most $d-3$ vertices. Consequently, the following hold.

$(a)$ Every vertex of $H_0$ has at most one neighbor in $Q^\circ$.

$(b)$ No two adjacent vertices of $H_0$ both have a neighbor in $Q^\circ$.

\p By Claim C, $G[Q^\circ]-E(C_1)$ is a forest each of whose components contains exactly one vertex of $C_1$. Let $a,b\in Q^\circ$, and let $a',b'$ be the vertices of $C_1$ in the components containing $a$ and $b$, respectively. The paths from $a$ to $a'$ and from $b'$ to $b$ in this forest, together with a shortest path from $a'$ to $b'$ in $C_1$, contain a path from $a$ to $b$ with at most $|Q^\circ\setminus V(C_1)|+\lfloor \tfrac{|V(C_1)|}{2} \rfloor+1=|Q^\circ|-\lceil \tfrac{|V(C_1)|}{2} \rceil+1\leq (h-1)-\lceil \tfrac{d}{2} \rceil+1=d-3$ vertices. If a vertex $w$ of $H_0$ has two neighbors $a,b\in Q^\circ$, then $w$ and such a path from $a$ to $b$ form a circuit of $G-v$ with at most $d-2$ vertices. If $ww'\in E(H_0)$, $a\in N(w)\cap Q^\circ$ and $b\in N(w')\cap Q^\circ$, then either $a=b$ and $ww'aw$ is a triangle, or $w$, $w'$ and such a path from $b$ to $a$ form a circuit of $G-v$ with at most $d-1$ vertices. Both contradict $g(G-v)=d\geq 4$. \q

\noindent\textbf{Claim E.} No vertex of $V(H_0)\setminus (Z\cup X_{H_0})$ has a neighbor in $Q$.

\p Suppose that $w\in V(H_0)\setminus (Z\cup X_{H_0})$ has a neighbor in $Q$. Since $w\notin N(v)$ and Claim D$(a)$, we have $N(w)\cap Q=\{a\}$ for some $a\in Q^\circ$. By Claim C, $G[Q]-E(C_1)$ is a forest each of whose components contains exactly one vertex of $C_1$. Hence at most $|Q\setminus V(C_1)|\leq h-d$ vertices of $C_1$ have degree at least three in $G[Q]$. Since $|V(C_1)|-(h-d)\geq 2d-h=d-\lfloor \frac{d-1}{2} \rfloor+2\geq 3$, there is a vertex $z\in V(C_1)\setminus\{a\}$ with $d_{G[Q]}(z)=2$. Since $v\notin V(C_1)$ and $v$ is adjacent to $v_3$, we have $z\notin\{v,v_3\}$, and thus $G[Q\setminus\{z\}]$ is a tree containing $v$ and $v_3$. So $T^*=G[(Q\setminus\{z\})\cup\{w\}]$ is an induced tree with $v,v_3\in V(T^*)$ and $V(T^*)\cap Z=\emptyset$. Since $\delta=3$ and Claim D$(a)$, $z$ has a neighbor in $H_0$, which is different from $w$. Since $H_0-w$ is connected, we have $G-V(T^*)=G[(V(H_0)\setminus\{w\})\cup\{z\}]$ is also connected. Hence $T^*$ is a good tree with $|V(T^*)|=|V(T)|+1$, which contradicts the maximality of $|V(T)|$. \q

Since $x_1$ is not a cut vertex of $H$, the graph $H_0=H-x_1$ is connected. By Claims A and E, applied with $U=Q$ and $J=H_0$, $H_0$ is a path with ends $v_1$ and $v_2$. Let $q$ be the neighbor of $v_1$ in $H_0$ (it is possible that $q=v_2$). Since $v_1$ has exactly one neighbor in $H_0$, $Q\setminus Q^\circ=\{v\}$ and $\delta=3$, the vertex $v_1$ has a neighbor in $Q^\circ$. Similarly, $q$ has a neighbor in $Q^\circ$. This contradicts Claim D$(b)$. Therefore $|V(T)|\geq h$, and this completes the proof of Lemma \ref{3.4}.
\q

\begin{lem}\label{3.5}
Let $G$ be a simple 3-connected graph with minimum degree $\delta=3$ and $v\in V(G)$ with $N(v)=\{v_1,v_2,v_3\}$. Suppose that $g(G-v)=d\geq 3$ and $h=d+\lfloor \frac{d-1}{2} \rfloor-2$.  If $G\ncong K_4$, then $G$ has a non-separating induced tree $T$ with $|V(T)|\geq h$ such that $v_3\in V(T)$ and $v,v_1,v_2\notin V(T)$.
\end{lem}

\p Since $g(G-v)=d$, every circuit of $G-v$ has length at least $d$. Since $\lfloor \frac{d-1}{2} \rfloor=\lceil \frac{d}{2} \rceil-1\leq \frac{d-1}{2}$, we have $2h\leq 3d-5$, and $h-\lceil \tfrac{d}{2} \rceil=d-3$.

Let $Z=\{v,v_1,v_2\}$ and $Z^\circ=\{v_1,v_2\}$. We call an induced tree $T$ of $G$ {\it good} if $v_3\in V(T)$, $V(T)\cap Z=\emptyset$ and $G-V(T)$ is connected. Since $G-v_3$ is connected, $G[\{v_3\}]$ is good, and we may choose a good tree $T$ with $|V(T)|$ as large as possible. Suppose, to the contrary, that $|V(T)|\leq h-1$. Then $|E(T)|\leq h-2$. Let $H=G-V(T)$. Since $N(v)=\{v_1,v_2,v_3\}$, the only neighbor of $v$ in $T$ is $v_3$, and no vertex other than $v_1,v_2,v_3$ is adjacent to $v$. For a graph $J$, let $X_J$ be the set of cut vertices of $J$.

\noindent\textbf{Claim A.} Let $U\subseteq V(G)\setminus Z$ with $v_3\in U$ such that $J=G-U$ is connected. If no vertex of $V(J)\setminus (Z\cup X_J)$ has a neighbor in $U$, then $V(J)=Z$.

\p Since $v_3\in U$, we have $N_J(v)=\{v_1,v_2\}$. Let $B$ be a block of $J$, and let $P_B$ be the set of vertices of $B$ which are not cut vertices of $J$. Every vertex of $P_B$ has all its neighbors in $V(B)\cup U$. Suppose that $S=P_B\setminus Z\neq\emptyset$. By hypothesis, no vertex of $S$ has a neighbor in $U$, and no vertex of $S$ is adjacent to $v$. So $N_G(S)\subseteq X_B\cup (Z^\circ\cap V(B))$. Since $N_G(S)$ separates $S$ from $U$ in $G$ and $G$ is 3-connected, we obtain the following.

$(\ast)$ If $P_B\not\subseteq Z$, then $|X_B\cup (Z^\circ\cap V(B))|\geq 3$.

If $X_J=\emptyset$, then $J$ is a block with $P_J=V(J)$, so $|X_B\cup (Z^\circ\cap V(B))|\leq 2$. By $(\ast)$, we get $V(J)=Z$, so the claim holds. Otherwise, $X_J\neq \emptyset$ and $J$ has at least two end blocks. Each end block $B$ contains exactly one cut vertex of $J$, and thus $P_B\neq\emptyset$. If one of end blocks satisfies $P_B\cap Z=\emptyset$, then $|X_B\cup (Z^\circ\cap V(B))|=1$, which contradicts $(\ast)$. Thus each end block satisfies $P_B\cap Z\neq\emptyset$. If $v$ is not a cut vertex of $J$, then the vertices $v,v_1,v_2$ lie in a common block of $J$ as $d_J(v)=2$. Thus at most one end block $B$ satisfies $P_B\cap Z\neq\emptyset$, a contradiction. Hence $v$ is a cut vertex of $J$, and the blocks of $J$ containing $v$ are the edges $vv_1$ and $vv_2$. Since each end block satisfies $P_B\cap Z\neq\emptyset$, we have $J$ has exactly two end blocks $B_1,B_2$, and we may assume that $P_{B_i}\cap Z=\{v_i\}$ for $i\in\{1,2\}$. This implies that $|X_{B_i}\cup (Z^\circ\cap V(B_i))|\leq 2$. By $(\ast)$, we get $P_{B_i}=\{v_i\}$ and $B_i$ is an edge $vv_i$. So $J$ is the path $v_1vv_2$, and thus $V(J)=Z$. 
\q

\noindent\textbf{Claim B.} $V(H)\neq Z$.

\p Suppose that $V(H)=Z$. Thus $V(G)=V(T)\cup Z$. If $|V(T)|=1$, then $|V(G)|=4$. Since $\delta=3$, we have $G\cong K_4$, a contradiction. Thus $|V(T)|\geq 2$. Since $|V(T)|\leq h-1$, we get $h\geq 3$. So $g(G-v)=d\geq 4$. Let $u\neq v_3$ be a leaf of $T$. Since $\delta=3$ and $u\notin N(v)$, the vertex $u$ is adjacent to both $v_1$ and $v_2$. Since $g(G-v)=d\geq 4$, we have $v_1v_2\notin E(G)$. Hence, $H$ is a path $v_1vv_2$. Since $\delta=3$, every neighbor of $v_i$ other than $v$ lies in $T$ for $i\in\{1,2\}$. Let $a_i\in (N(v_i)\cap V(T))\setminus\{u\}$. Then $v_1uTa_1v_1$, $v_2uTa_2v_2$ and $v_1a_1Ta_2v_2uv_1$ are circuits of $G-v$ (where $a_1Ta_2$ is trivial if $a_1=a_2$) of lengths $|E(uTa_1)|+2$, $|E(uTa_2)|+2$ and $|E(a_1Ta_2)|+4$, respectively. Since every edge of the tree $T$ lies on at most two of the paths $uTa_1$, $uTa_2$ and $a_1Ta_2$, the sum of these three lengths is at most $2|E(T)|+8\leq 2h+4\leq 3d-1$. This contradicts $g(G-v)=d\geq 4$. 
\q

By Claims A and B, applied with $U=V(T)$ and $J=H$, there is a vertex $x_1\in V(H)\setminus (Z\cup X_H)$ with a neighbor in $T$. Let $Q=V(T)\cup\{x_1\}$, and $H_0=H-x_1$.

\noindent\textbf{Claim C.} $x_1$ has exactly two neighbors $y_1,z_1$ in $T$. Consequently, $|Q|\leq h$, the graph $G[Q]$ is connected and has a unique circuit $C_1=x_1y_1Tz_1x_1$.

\p If $x_1$ has exactly one neighbor in $T$, then $G[V(T)\cup\{x_1\}]$ is a good tree as $H-x_1$ is connected, which contradicts the maximality of $|V(T)|$. If $x_1$ has three neighbors $y,z,w$ in $T$, then $y,z,w\in V(T)$, and $x_1yTzx_1$, $x_1zTwx_1$ and $x_1yTwx_1$ are circuits of $G-v$ of lengths $|E(yTz)|+2$, $|E(zTw)|+2$ and $|E(yTw)|+2$, respectively. Since every edge of the tree $T$ lies on at most two of the paths $yTz$, $zTw$ and $yTw$, the sum of their lengths is at most $2|E(T)|+6\leq 2h+2\leq 3d-3$, which contradicts the fact that each of them has length at least $d$. Hence $x_1$ has exactly two neighbors $y_1,z_1$ in $T$, and $y_1,z_1\in V(T)$. Since $T$ is a tree and $H-x_1$ is connected, we have $|Q|=|V(T)|+1\leq h$, $G[Q]$ is connected and has a unique circuit $C_1=x_1y_1Tz_1x_1$. \q

\noindent\textbf{Claim D.} Every vertex of $V(G)\setminus(Q\cup\{v\})$ has at most one neighbor in $Q$.

\p By Claim C, $G[Q]-E(C_1)$ is a forest each of whose components contains exactly one vertex of $C_1$. Let $a,b\in Q$, and let $a',b'$ be the vertices of $C_1$ in the components containing $a$ and $b$, respectively. The paths from $a$ to $a'$ and from $b'$ to $b$ in this forest, together with a shortest path from $a'$ to $b'$ in $C_1$, contain a path from $a$ to $b$ with at most $|Q\setminus V(C_1)|+\lfloor \tfrac{|V(C_1)|}{2} \rfloor+1=|Q|-\lceil \tfrac{|V(C_1)|}{2} \rceil+1\leq h-\lceil \tfrac{d}{2} \rceil+1=d-2$ vertices. If a vertex $w\in V(G)\setminus(Q\cup\{v\})$ had two neighbors $a,b\in Q$, then $w$ and such a path from $a$ to $b$ would form a circuit of $G-v$ with at most $d-1$ vertices, which contradicts $g(G-v)=d$. \q

\noindent\textbf{Claim E.} No vertex of $V(H_0)\setminus (Z\cup X_{H_0})$ has a neighbor in $Q$.

\p Suppose that $w\in V(H_0)\setminus (Z\cup X_{H_0})$ has a neighbor in $Q$. By Claim D, $N(w)\cap Q=\{a\}$ for some $a\in Q$. By Claim C, $G[Q]-E(C_1)$ is a forest each of whose components contains exactly one vertex of $C_1$. Hence at most $|Q\setminus V(C_1)|\leq h-d$ vertices of $C_1$ have degree at least three in $G[Q]$. Since $|V(C_1)|-(h-d)\geq 2d-h=d-\lfloor \frac{d-1}{2} \rfloor+2\geq 3$, there is a vertex $z\in V(C_1)\setminus\{a,v_3\}$ with $d_{G[Q]}(z)=2$. Thus $G[Q\setminus\{z\}]$ is a tree containing $v_3$ and $T^*=G[(Q\setminus\{z\})\cup\{w\}]$ is an induced tree with $v_3\in V(T^*)$ and $V(T^*)\cap Z=\emptyset$. Since $\delta=3$ and Claim D, $z$ has a neighbor in $H_0$, which is different from $w$. Since $H_0-w$ is connected, we have $G-V(T^*)=G[(V(H_0)\setminus\{w\})\cup\{z\}]$ is also connected. Hence $T^*$ is a good tree with $|V(T^*)|=|V(T)|+1$, which contradicts the maximality of $|V(T)|$. \q

Since $x_1$ is not a cut vertex of $H$, the graph $H_0=H-x_1$ is connected. By Claims A and E, applied with $U=Q$ and $J=H_0$, we have $V(H_0)=Z$, and thus $V(G)=Q\cup Z$. By Claim C, $G[Q]$ has exactly $|Q|$ vertices and $|Q|$ edges. Since $\delta=3$, there are at least $3|Q|-2|Q|=|Q|$ edges joining $Q$ to $Z$. On the other hand, by Claim D, each of $v_1,v_2$ has at most one neighbor in $Q$. Since $v_3$ is the only neighbor of $v$ in $Q$ and $V(H_0)=Z$, we have $|Q|\leq 3$. Since $|V(G[Q])|=|Q|\geq |V(C_1)|\geq d\geq 3$, we get $|Q|=d=3$. Thus $h=2$. But $|Q|=3>2=h$, which contradicts Claim C. Therefore $|V(T)|\geq h$, and this completes the proof of Lemma \ref{3.5}.
\q

The prism graph is the graph obtained from two vertex-disjoint
triangles by adding a perfect matching between their vertex sets.

\begin{theo}\label{3.6}
Let $G$ be a simple 3-connected graph and $T=\{e,f,g\}$ be a 3-cocircuit in $G$. Suppose that $g(G-T)\geq d\geq 3$ and $h=d+\lfloor \frac{d-1}{2} \rfloor-2$. If $G \ncong K_4$, then there exist bonds $C_i$ and $D_j$ of $G$ for $i\in \{1,\ldots, h\}$, and $j\in \{1,\ldots, h-1\}$ such that $C_i\cap T=\{e\}$ and $D_j\cap T=\{f,g\}$, with $|C_1|,\ldots,|C_{h}|$ pairwise distinct and $|D_1|,\ldots,|D_{h-1}|$ pairwise distinct.
\end{theo}

\p Let $F$ be a component of $G-T$, where we choose $F$ to be a single vertex if $T$ is the set of edges incident with a vertex. Let $G_1$ be obtained from $G$ by contracting $F$ to a single vertex $v$. Let $v_1$, $v_2$ and $v_3$ be the vertices of $G_1$ which are joined to $v$ by $f,g$ and $e$, respectively. Then $G_1$ is also simple and 3-connected, $T=\{e,f,g\}$ is a 3-cocircuit in $G_1$ and $g(G_1-T)\geq d\geq 3$. Furthermore, every cocircuit of $G_1$ is also a cocircuit of $G$. If $F$ is a vertex, then $G_1=G$. If $F$ has more than one vertex, then both components have at least three vertices, and thus 3-connectivity implies that three edges of $T$ have distinct ends in each component. Then a 2-cut of $G_1$ avoiding $v$ would separate $G$. If $G_1-\{v,x\}$ were disconnected, one component would meet at most one edge of $T$, implying that there would be a vertex cut of $G$ of size at most two. For either choice of the contracted component used below, put $d'=g(G_1-v)$. Since $G_1-v$ is a component of $G-T$, we have $d'\geq d$, and hence $h'=d'+\left\lfloor\frac{d'-1}{2}\right\rfloor-2 \geq d+\left\lfloor\frac{d-1}{2}\right\rfloor-2=h.$ We apply Lemmas~\ref{3.4} and~\ref{3.5} with the actual girth $d'$ and retain only the required number of bonds. By Lemma \ref{3.4} and the proof of Lemma \ref{2.7}, where the deleted leaves are chosen outside $R=\{v,v_3\}$, and retaining $h-1\leq h'-1$ bonds, there exist bonds $D_j$ of $G$ for $j\in \{1,\ldots, h-1\}$ such that $D_j\cap T=\{f,g\}$, with  $|D_1|,\ldots,|D_{h-1}|$ pairwise distinct. If $G_1\ncong K_4$, then, combining Lemma \ref{3.5} with the proof of Lemma \ref{2.7} with $R=\{v_3\}$, there exist bonds $C_i$ of $G$ for $i\in \{1,\ldots, h\}$ such that $C_i\cap T=\{e\}$, with $|C_1|,\ldots,|C_{h}|$ pairwise distinct. If $G_1\cong K_4$, we contract the other component of $G-T$ instead. This fails only if both components of $G-T$ are triangles, that is, $G$ is the prism and $T$ is its perfect matching. In this case, $d=3$ and $h=2$. Writing the two triangles as $a_1a_2a_3$ and $b_1b_2b_3$ with $e=a_3b_3$, the bonds $\partial(\{a_3\})$ and $\partial(\{a_1,a_3,b_1\})$ have sizes $3$ and $5$ and meet $T$ exactly in $\{e\}$.
\q

Our main Theorem \ref{1.2} follows immediately from Theorems \ref{3.3} and \ref{3.6}.

\section{Binary matroids}

We will use the decomposition theorem of Seymour~\cite{Se1980} to extend Theorems \ref{1.1} and \ref{1.2} to binary matroids which do not have both an $F_7$-minor and an $F_7^*$-minor. We first need some lemmas to construct circuits in 2-sums and 3-sums of binary matroids from Hochst\"{a}ttler and Jackson~\cite{Hoch1998}.

\begin{lem}\label{4.1} {\rm (Hochst\"{a}ttler and Jackson~\cite{Hoch1998})}
Suppose that $M$ is a binary matroid and $M=M_1\oplus_3 M_2$ for minors $M_1$ and $M_2$ of $M$, where $E(M_1)\cap E(M_2)=X$ for some triangle $X=\{x,y,z\}$ of $M_i$, $1\leq i\leq 2$. Let $C_i$ be a circuit of $M_i$ such that $X_1=C_1\cap X=\{x,y\}$ and $X_2=C_2\cap X=\{z\}$. Then $C_1\Delta C_2\Delta X$ is a circuit of $M$.
\end{lem}

\begin{lem}\label{4.2} {\rm (Hochst\"{a}ttler and Jackson~\cite{Hoch1998})}
Let $M$ be a connected binary matroid. Suppose that $e\in E(M)$ and $X=\{x,y,z\}$ is a triangle in $M-e$. Then $M$ has a circuit containing $e$ and exactly one element of $X$.
\end{lem}

\begin{lem}\label{4.3} {\rm (Hochst\"{a}ttler and Jackson~\cite{Hoch1998})}
Let $M$ be a vertically 3-connected binary matroid. Suppose that $T=\{e,f,g\}$ and $X=\{x,y,z\}$ are two disjoint triangles in $M$. Then $M$ has a circuit $C$ with $C\cap T=\{f,g\}$ and $|C\cap X|=1$ and a circuit $D$ with $D\cap T=\{e\}$ and $|D\cap X|=2$.
\end{lem}

\begin{transferlem}
Let $M=N_1\oplus_k N_2$ be a 2-sum or 3-sum of binary matroids with common set $S$. If $M_i$ is obtained from $N_i$ by deleting elements parallel to members of $S$ retaining $S$, then every cocircuit of $M_i$ avoiding $S$ is a cocircuit of $M$.
\end{transferlem}
\p Let $P_i$ be the deleted set. Let $Y$ be a cocircuit of $M_i$ avoiding $S$. Then there exists a cocircuit $D$ of $N_i$ such that $D\setminus P_i=Y$. Since $D$ avoids $S$, orthogonality with each pair $\{p,s\}$ implies that $D=Y$, where $p\in P_i$ is parallel to $s\in S$. Since $S$ contains no cocircuit of $N_j$, we have $N_i/S=M/(E(N_j)\setminus S)$ for $j\neq i$. Since $Y$ is a cocircuit of $N_i/S$, it is a cocircuit of $M$. \q

Now, we will prove a linear lower bound of the size of the circuit spectra for simple connected binary matroids not having both an $F_7$-minor and an $F_7^*$-minor, and having minimum cocircuit size $d\geq 3$. Our first main theorem, Theorem \ref{1.3}, follows from the next theorem. We use the method of \cite{Hoch1998}, which uses the Decomposition theorem of Seymour \cite{Se1980}.

\begin{theo}\label{4.4} 
Let $M$ be a simple binary matroid that does not have both an $F_7$-minor and an $F_7^*$-minor.

$(a)$ Suppose $M$ is connected, $|E(M)|\geq 2$, $e\in E(M)$ and every cocircuit $Y$ of $M$ with $e\notin Y$ has size at least $d\geq 3$. Then the size of the circuit $e$-spectrum of $M$ is at least $d-2$.

$(b)$ Suppose $M$ is 3-connected, $|E(M)|\geq 4$, $T=\{e,f,g\}$ is a triangle in $M$ and every cocircuit $Y$ of $M$ with $T\cap Y=\emptyset$ has size at least $d\geq 3$. If $M\ncong F_7$ when $d=4$, then there exist circuits $C_i,D_i$ of $M$ for $i\in \{1,\ldots, d-2\}$ such that $C_i\cap T=\{e\}$ and $D_i\cap T=\{f,g\}$, with $|C_1|,\ldots,|C_{d-2}|$ pairwise distinct and $|D_1|,\ldots,|D_{d-2}|$ pairwise distinct.
\end{theo}

\p Suppose, to the contrary, that $M$ is a counterexample chosen such that $|E(M)|$ is as small as possible. Note that the hypotheses $|E(M)|\geq 2$ in $(a)$ and $|E(M)|\geq 4$ in $(b)$ cannot be omitted. 

\vspace{2mm}
{\bf Claim 1} {\em $M$ is 3-connected.}
\vspace{2mm}

\p Suppose that $M$ is not 3-connected. By Theorem \ref{2.3}, $M=N_1\oplus_2 N_2$ for minors $N_1$ and $N_2$ of $M$, where $E(N_1)\cap E(N_2)=\{x\}$ and $|E(N_i)|\geq 3$ for $i\in \{1,2\}$. Let $M_i$ be the simple matroid obtained from $N_i$ by deleting any elements that are parallel to $x$ ($i\in \{1,2\}$). It follows that $M_i$ and $N_i$ are both minors of $M$. Since $M$ is connected, both $M_1$ and $M_2$ are connected. Suppose that $(a)$ does not hold for some $e\in E(M)$. By symmetry, assume that $e\in E(N_2)$. Since every cocircuit of $M$ which avoids $e$ has size at least $d$, by the cocircuit-transfer lemma, every cocircuit of $M_1$ which avoids $x$ has size at least $d$. Note that $|E(M_1)|<|E(M)|$. Moreover, $|E(M_1)|\geq 2$. In fact, two elements of $N_1$ parallel to $x$ would form a 2-element circuit of $N_1$ avoiding $x$, and hence a 2-element cycle of $M$, contrary to $M$ being simple. So $|E(M_1)|\geq|E(N_1)|-1\geq 2$. By the minimality of $|E(M)|$ and $(a)$, the size of the circuit $x$-spectrum of $M_1$ is at least $d-2$. Suppose that $C_i\in \mathcal{C}(M_1)$ contains $x$ for $i\in \{1,\ldots d-2\}$, with the $|C_i|$ pairwise distinct. Since $N_2$ is connected, $N_2$ has a circuit $Q$ containing both $x$ and $e$. Thus $C_1\Delta Q, \ldots, C_{d-2}\Delta Q$ are $d-2$ circuits of $M$ with distinct lengths containing $e$. This implies that the size of the circuit $e$-spectrum of $M$ is at least $d-2$. This contradicts the fact that $M$ does not satisfy $(a)$. 
\q

\vspace{2mm}
{\bf Claim 2} {\em $M$ is regular.}
\vspace{2mm}

\p
By Theorem \ref{2.5} and Claim 1, it follows that $M$ is regular, or is isomorphic to $F_7$ or $F^*_7$. Note that $F_7$ has only 3-element and 4-element circuits, and $F^*_7$ has only 4-element circuits and it does not have any triangles. It is easily checked that $(a)$ holds for $F_7$ and $F_7^*$. Since $F_7^*$ has no triangles, we only need to consider $F_7$ in $(b)$. In this case, $d\leq4$, and $d=4$ is excluded by hypothesis. Hence $d=3$, and the required circuits are easily found. This contradicts the fact that $M$ does not satisfy $(a)$ or $(b)$. Thus $M$ is regular.
\q

\vspace{2mm}
{\bf Claim 3} {\em $M$ cannot be expressed as a $3$-sum.}
\vspace{2mm}

\p
First suppose that $M=N_1\oplus_3 N_2$ for minors $N_1$ and $N_2$ of $M$, where $E(N_1)\cap E(N_2)=X$ for some triangle $X=\{x,y,z\}$ of $N_i$ for $i\in \{1,2\}$ and $|E(N_i)|\geq 7$. Let $M_i$ be the simple matroid obtained from $N_i$ by deleting any elements that are parallel to $x,y,z$. By Claim 1, $M$ is 3-connected. It follows from \cite[(4.3)]{Se1980} that $N_i$ is vertically 3-connected and $M_i$ is 3-connected for $i\in \{1,2\}$. By Claim 2, $M$ is regular. Note that $N_i$ and $M_i$ are both minors of $M$, Thus $N_i$ and $M_i$ are both regular for $i\in \{1,2\}$. Next, we consider the following cases.

Suppose that $(a)$ does not hold for some $e\in E(M)$. By symmetry, we may assume that $e\in E(N_2)$. If $M=M(K_4)$, then it is easily checked that $(a)$ holds, a contradiction. Suppose that $M\ncong M(K_4)$. By Lemma \ref{4.2}, $N_2$ has a circuit $Q$ containing $e$ and exactly one element, say $z$, of $X$. Since every cocircuit of $M$ which avoids $e$ has size at least $d$, by the cocircuit-transfer lemma, every cocircuit of $M_1$ which avoids $X$ has size at least $d$. Note that $|E(M_1)|<|E(M)|$. Moreover, $|E(M_1)|\geq 4$ and $M_1\ncong F_7$. In fact, two elements of $E(N_1)-X$ parallel in $N_1$ to the same element of $X$ would form a 2-element circuit of $N_1$ avoiding $X$, and hence a 2-element cycle of $M$, contrary to $M$ being simple. So $|E(M_1)|\geq|E(N_1)|-3\geq 4$. Also, $M_1$ is regular, so $M_1\ncong F_7$. The same holds for $M_2$. By the minimality of $|E(M)|$ and $(b)$, there exists $C_i\in \mathcal{C}(M_1)$ for $i\in \{1,\ldots d-2\}$ such that $C_i\cap X=\{x,y\}$, with the $|C_i|$ pairwise distinct. By Lemma \ref{4.1}, $C_1\Delta Q\Delta X,\ldots, C_{d-2}\Delta Q\Delta X$ are $d-2$ circuits of $M$ with distinct lengths containing $e$. This contradicts the fact that $M$ does not satisfy $(a)$.

Now suppose that $(b)$ does not hold for some triangle $T=\{e,f,g\}$ of $M$. We may assume that $N_1$ and $N_2$ have been chosen such that $T\cap E(N_2)$ is as large as possible.

We first show that $T\subseteq E(N_2)$. Suppose that $T\cap E(N_1)\neq \emptyset$. By the choice of $N_1$ and $N_2$, we have $|T\cap E(N_1)|=1$ and $|T\cap E(N_2)|=2$. Let $S_i=E(N_i)-X$ and $T\cap E(N_1)=\{p\}$. Since $M=N_1\oplus_3 N_2$ and $M$ is 3-connected, we have $r(S_1)+r(S_2)=r(M)+2$. Let $S'_1=S_1-\{p\}$ and $S'_2=S_2\cup \{p\}$. Since  $r(S'_1)\leq r(S_1)$ and $p\in\operatorname{cl}_M(S_2)$, we have $r(S'_1)\leq r(S_1)$ and $r(S'_2)=r(S_2)$. Thus $r(S'_1)+r(S'_2)\leq r(M)+2$. As $M$ is 3-connected, it must be true that $r(S'_1)+r(S'_2)=r(M)+2$. If $|S_1|\geq 5$, then $|S'_1|\geq 4$. Since $(S'_1,S'_2)$ is an exact 3-separation of $M$, it follows from Theorem \ref{2.4} that there must exist two matroids $N'_1$ and $N'_2$, both of which are minors of $M$,  such that $M=N'_1\oplus_3 N'_2$, where $E(N'_i)\cap E(M)=S'_i$. This contradicts the maximality of $T\cap E(N_2)$ and hence $|S_1|=4$ as $|E(N_i)|\geq 7$. Since $M$ is simple and binary, we have $r(S_1)\geq 3$. Thus $r(S'_2)=r(S_2)\leq r(M)-1$. Hence, $S'_2$ is contained in a hyperplane of $M$, and $S'_1$ contains a cocircuit of $M$. Since $|S'_1|=3$ and $S'_1\cap T=\emptyset$, we have $d=3$. 

Since $d=3$, it suffices to construct the two required circuits.
Since $M$ is 3-connected, $M\setminus e$ is connected. So it contains a circuit $D$ through $f$ and $g$. Put $C=D\Delta T$.
Since $M$ is binary, $C$ is a nonempty cycle. Since any circuit $Q\subseteq C$ must contain $e$, $Q\Delta T$ is a nonempty cycle contained in $D$. So $Q\Delta T=D$ and $Q=C$. Thus $C$ is a circuit with $C\cap T=\{e\}$ and $D\cap T=\{f,g\}$. This proves $(b)$, a contradiction. Therefore, $T\subseteq E(N_2)$.

Since $T\subseteq E(N_2)\setminus X$, $T$ is also a triangle of $N_2$. Applying Lemma \ref{4.3} to $N_2$, we deduce that $N_2$ has a circuit $Q'$ with $Q'\cap T=\{e\}$ containing exactly two elements, say $a,b$, of $X$ and has a circuit $Q''$ with $Q''\cap T=\{f,g\}$ containing exactly one element, say $c$, of $X$. It is possible that $c\in \{a,b\}$. Since every cocircuit of $M$ which avoids $T$ has size at least $d$, by the cocircuit-transfer lemma, every cocircuit of $M_1$ which avoids $X$ has size at least $d$. Applying $(b)$ to $M_1$ (recall that $|E(M_1)|\geq 4$ and $M_1\ncong F_7$), we deduce that there exists $C'_i,C''_i\in \mathcal{C}(M_1)$ for $i\in \{1,\ldots, d-2\}$ such that $C'_i\cap X=X\setminus \{a,b\}$ and $C''_i\cap X=X\setminus \{c\}$, with the $|C'_i|$ pairwise distinct and the $|C''_i|$ pairwise distinct. By Lemma \ref{4.1}, $C'_1\Delta Q'\Delta X,\ldots, C'_{d-2}\Delta Q'\Delta X$ are $d-2$ circuits of $M$ with distinct lengths containing $e$ and $C''_1\Delta Q''\Delta X,\ldots, C''_{d-2}\Delta Q''\Delta X$ are $d-2$ circuits of $M$ with distinct lengths containing $f$ and $g$. This contradicts the fact that $M$ does not satisfy $(b)$.
\q

By Claims 1, 2, 3 and Theorem \ref{2.6}, $M$ is graphic or cographic or is isomorphic to $R_{10}$. If $M=R_{10}$, then $d\leq 4$, and every element belongs to a 4-element circuit and a 6-element circuit, and it has no triangles. Thus both $(a)$ and $(b)$ hold. If $M$ is graphic or cographic, then it follows from Corollaries~\ref{3.1} and~\ref{3.2} (graphic case, since every vertex star of $G$ is a cocircuit of $M$) and Theorem~\ref{1.2} (cographic case) that both $(a)$ and $(b)$ hold (if $M\cong M(K_4)$, we use that $M$ is graphic and apply Corollary~\ref{3.2}). In the cographic case of $(a)$, we apply Theorem~\ref{1.2}$(a)$ with $d$ replaced by $g(G-e)\ge d$. Here we use $|E(M)|\geq 2$ in $(a)$ and $|E(M)|\geq 4$ in $(b)$, which exclude $U_{1,1}\cong M(K_2)$ and $U_{2,3}\cong M(K_3)$, respectively. Corollaries~\ref{3.1}, \ref{3.2} and Theorem~\ref{1.2} do not apply to them, as $K_2$ is not 2-connected, $K_3$ is not 3-connected, and $U_{1,1}^*$ and $U_{2,3}^*$ are the cycle matroids of a loop and of three parallel edges. This contradicts the fact that $M$ does not satisfy $(a)$ or $(b)$. This completes the proof of Theorem \ref{4.4}.

\begin{theo}\label{4.5}
Let $M$ be a nonempty simple connected binary matroid that does not have both an $F_7$-minor and an $F_7^*$-minor. Suppose that every cocircuit of $M$ has size at least $d\geq 3$. Then the size of the circuit spectrum of $M$ is at least $d-1$, unless $M\cong F_7,F^*_7,R_{10}$. In the exceptional cases, the size of the circuit spectrum of $M$ is exactly $d-2$ if $d$ is the minimum cocircuit size of $M$.
\end{theo}

\p For $M\cong F_7,F_7^*$, or $R_{10}$, the assertion follows by direct inspection. Their circuit spectra are $\{3,4\}$, $\{4\}$, and $\{4,6\}$, respectively, and their minimum cocircuit sizes are $4$, $3$, and $4$. If $d$ equals the minimum cocircuit size, then $|C(M)|=d-2$. If $d\geq 3$ is a smaller admissible integer, then $d=3$ and $M\cong F_7$ or $R_{10}$, so $|C(M)|=2=d-1$. First, we will show that when $d=3$, we get the size of the circuit spectrum of $M$ at least $d-1$ unless $M\cong F^*_7$. Note that $M$ is simple and cosimple ($d\geq 3$). If $M$ has the circuit spectrum of size one, then by Theorem \ref{Murty}, $M$ is the dual of projective geometry $PG(r,2)$ $(r\geq 2)$ or the dual of the Affine geometry $AG(r+1,2)$ $(r\geq 2)$. So if $M\ncong F^*_7$, then $M\cong PG(r,2)^*$ $(r\geq 3)$ or $M\cong AG(r+1,2)^*$ $(r\geq 2)$. It follows that $M$ contains both $F_7$ and its dual as a minor, which contradicts the assumption. Thus $M\cong F^*_7$. 

Now, it suffices to show that Theorem \ref{4.5} holds for $M\ncong F_7,F^*_7,R_{10}$ and $d\geq 4$. We proceed by contradiction. Suppose that the theorem is false.

\vspace{2mm}
{\bf Claim 1} {\em $M$ is 3-connected and regular}.
\vspace{2mm}

\p Suppose that $M$ is not 3-connected. By Theorem \ref{2.3}, $M=N_1\oplus_2 N_2$ for minors $N_1$ and $N_2$ of $M$, where $E(N_1)\cap E(N_2)=\{x\}$ and $|E(N_i)|\geq 3$ for $i\in \{1,2\}$. Let $M_i$ be the simple matroid obtained from $N_i$ by deleting any elements which are parallel to $x$. It follows that $M_i$ and $N_i$ are both minors of $M$. Since $M$ is connected, $M_i$ is connected for $i\in \{1,2\}$. Since every cocircuit of $M$ has size at least $d$, by the cocircuit-transfer lemma, every cocircuit of $M_i$ which avoids $x$ has size at least $d$. Since at most one element parallel to $x$ is deleted, we have $|E(M_i)|\geq|E(N_i)|-1\geq2$. Applying Theorem \ref{4.4} $(a)$ to $M_i$, we deduce that the size of the circuit $x$-spectrum of $M_i$ is at least $d-2$. Suppose that $C_i\in \mathcal{C}(M_1)$ and $D_i\in \mathcal{C}(M_2)$ contain $x$ for $i\in \{1,\ldots d-2\}$, with the $|C_i|$ pairwise distinct and the $|D_i|$ pairwise distinct. Thus $C_i\Delta D_j$ is a circuit of $M$ for $i,j\in \{1,\ldots,d-2\}$. Let $A$ and $B$ be the set of $|C_i-\{x\}|$ and $|D_i-\{x\}|$ for $i\in \{1,\ldots d-2\}$, respectively. By Lemma~\ref{Cauchy}, the size of the circuit spectrum of $M$ is at least $|A|+|B|-1\geq d-2+d-2-1\geq d-1$ for $d\geq 4$. This contradicts the fact that $M$ does not satisfy Theorem \ref{4.5}. Thus $M$ is 3-connected. By Theorem \ref{2.5}, it follows that $M$ is regular, or is isomorphic to $F_7$ or $F^*_7$. Note that $M\ncong F_7,F^*_7$. Thus $M$ is regular.
\q

\vspace{2mm}
{\bf Claim 2} {\em $M$ cannot be expressed as a $3$-sum.}
\vspace{2mm}

\p
Suppose that $M=N_1\oplus_3 N_2$ for minors $N_1$ and $N_2$ of $M$, where $E(N_1)\cap E(N_2)=X$ for some triangle $X=\{x,y,z\}$ of $N_i$ for $i\in \{1,2\}$ and $|E(N_i)|\geq 7$. Let $M_i$ be the simple matroid obtained from $N_i$ by deleting any elements which are parallel to $x,y,z$. By Claim 1, $M$ is 3-connected and regular. By \cite[(4.3)]{Se1980}, $N_i$ is vertically 3-connected and $M_i$ is 3-connected for $i\in \{1,2\}$. Since $N_i$ and $M_i$ are minors of $M$, $N_i$ and $M_i$ are regular. Since every cocircuit of $M$ has size at least $d$, by the cocircuit-transfer lemma, every cocircuit of $M_i$ which avoids $X$ has size at least $d$. The same parallel-element argument gives $|E(M_i)|\geq|E(N_i)|-3\geq4$, and regularity gives $M_i\ncong F_7$. Applying Theorem \ref{4.4} $(b)$ to $M_i$, we deduce that there exist $C_i\in \mathcal{C}(M_1)$ and $D_i\in \mathcal{C}(M_2)$ for $i\in \{1,\ldots d-2\}$ such that $C_i\cap X=\{z\}$ and $D_i\cap X=\{x,y\}$, with the $|C_i|$ pairwise distinct and the $|D_i|$ pairwise distinct. By Lemma \ref{4.1}, $C_i\Delta D_j\Delta X$ is a circuit of $M$ for $i,j\in \{1,\ldots,d-2\}$. Let $A$ and $B$ be the set of $|C_i-\{z\}|$ and $|D_i-\{x,y\}|$ for $i\in \{1,\ldots d-2\}$, respectively. By Lemma~\ref{Cauchy}, the size of the circuit spectrum of $M$ is at least $|A|+|B|-1\geq d-2+d-2-1\geq d-1$ for $d\geq 4$. This contradicts the fact that $M$ does not satisfy Theorem \ref{4.5}.

By Claims 1, 2 and Theorem \ref{2.6}, $M$ is graphic or cographic or is isomorphic to $R_{10}$. Note that $M\ncong R_{10}$. If $M$ is graphic or cographic, then it follows from Corollary~\ref{3.1} (graphic case) and Theorem~\ref{1.2}$(a)$ with $d$ replaced by $g(G-e)\ge d$ for an edge $e$ of $G$ (cographic case) that the size of the circuit spectrum of $M$ is at least $d-1$. This contradicts the fact that $M$ does not satisfy Theorem \ref{4.5}. This completes the proof of Theorem \ref{4.5}. \q

\section{Acknowledgments}

H. Wu’s research was supported by the 2025 College of Liberal Arts Senior Faculty Seed Funding from the University
of Mississippi.

\section*{Declaration}

All ideas and proof strategies were developed by the authors without using any AI tools. This work is part of the work in the third author’s dissertation, which did not use any AI tools. Claude was used for grammar checking, proofreading, and improving the clarity of the exposition. The authors take full responsibility for the accuracy and originality of the paper.

\end{document}